\documentclass[12pt,a4paper]{article}
\usepackage{}
\usepackage{amsfonts}
\usepackage{amssymb}
\usepackage{mathrsfs}
\usepackage{amsmath}
\usepackage{amsthm}
\usepackage{enumerate}
\usepackage{cite}
\usepackage[all]{xy}
\allowdisplaybreaks
\newtheorem{defn}{Definition}[section]
\newtheorem{thm}[defn]{Theorem}
\newtheorem{lem}[defn]{Lemma}
\newtheorem{prop}[defn]{Proposition}
\newtheorem{cor}[defn]{Corollary}
\newtheorem{eg}[defn]{Example}
\newtheorem{re}[defn]{Remark}

\newcommand\relphantom[1]{\mathrel{\phantom{#1}}}

\newcommand{\bdefn}{\begin{defn}}
\newcommand{\edefn}{\end{defn}}
\newcommand{\bthm}{\begin{thm}}
\newcommand{\ethm}{\end{thm}}
\newcommand{\blem}{\begin{lem}}
\newcommand{\elem}{\end{lem}}
\newcommand{\bprop}{\begin{prop}}
\newcommand{\eprop}{\end{prop}}
\newcommand{\bcor}{\begin{cor}}
\newcommand{\ecor}{\end{cor}}
\newcommand{\beg}{\begin{eg}}
\newcommand{\eeg}{\end{eg}}
\newcommand{\bre}{\begin{re}}
\newcommand{\ere}{\end{re}}
\newcommand{\bpf}{\begin{proof}}
\newcommand{\epf}{\end{proof}}

\newcommand{\Der}{{\rm Der}}

\newcommand{\id}{{\rm id}}

\newcommand{\K}{\mathbb{K}}

\newcommand{\g}{\mathfrak{g}}
\newcommand{\A}{\mathscr{A}}

\newcommand{\x}{\mathscr{X}}

\newcommand{\benu}{\begin{enumerate}}
\newcommand{\eenu}{\end{enumerate}}
\newcommand{\bc}{\begin{center}}
\newcommand{\ec}{\end{center}}
\newcommand{\bea}{\begin{eqnarray}}
\newcommand{\eea}{\end{eqnarray}}
\newcommand{\Bea}{\begin{eqnarray*}}
\newcommand{\Eea}{\end{eqnarray*}}
\newcommand{\beq}{\begin{equation}}
\newcommand{\eeq}{\end{equation}}
\newcommand{\Beq}{\begin{equation*}}
\newcommand{\Eeq}{\end{equation*}}
\newcommand{\bspl}{\begin{split}}
\newcommand{\espl}{\end{split}}

\numberwithin{equation}{section}

\begin{document}
\title{{\bf  $n$-ary Hom-Nambu algebras}}
\author{ Jun Zhao,  Liangyun Chen
 \date{{\small {School of Mathematics and Statistics, Northeast Normal
 University,\\
Changchun 130024, China}\\}}}

\maketitle
\date{}

\begin{abstract}

 In this paper, we define $\omega$-derivations, and study some properties of $\omega$-derivations, with its properties we can structure a new $n$-ary Hom-Nambu algebra from an $n$-ary Hom-Nambu algebra. In addition, we also give derivations and representations of $n$-ary Hom-Nambu algebras.
\bigskip

\end{abstract}
\renewcommand{\thefootnote}{\fnsymbol{footnote}}
\footnote[0]{ Corresponding author(L. Chen): chenly640@nenu.edu.cn.}
\footnote[0]{Supported by NNSF of China (Nos. 11171055 and
11471090). }

\section{Introduction}
Leibniz $n$-algebras were introduced by Casas, J. M. in [2], $n$-ary Hom-Nambu algebra is a generalization of the Leibniz $n$-algebra. An $n$-ary Hom-Nambu algebra is a triple $(\g, [\cdot,\cdots,\cdot],\alpha)$ consisting of a vector space $\g$, a multilinear map $[\cdot,\cdots,\cdot]: \underbrace{\g\times\cdots\times \g}_{n}\rightarrow \g$ and a family $\alpha=(\alpha_{i})_{1\leq i\leq n-1}$ of linear maps $\alpha_{i}: \g\rightarrow \g,$ satisfying
\begin{align*}
[[x_1,\cdots,x_n],\alpha_{1}(y_1),\cdots,\alpha_{n-1}(y_{n-1})]\\
=\sum_{i=1}^n[\alpha_{1}(x_1),\cdots,\alpha_{i-1}(x_{i-1}),&[x_i,y_1,\cdots, y_{n-1}],\alpha_{i}(x_{i+1}),\cdots,\alpha_{n-1}(x_n)].
\end{align*}
When $\alpha_{i}=\id$, it becomes a Leibniz $n$-algebra. For $n=2$, one recovers Hom-Leibniz algebras, the specific content about Hom-Leibniz algebras can be seen in [1,3].

In this paper, the main result is to study $\omega$-derivations of a vector space, and how to structure a new $n$-ary Hom-Nambu algebra. In addition, we also introduce the derivations and representations similar with $n$-ary multiplicative Hom-Nambu-Lie superalgebras and $n$-Lie superalgebras in [4,6] as its application.

The paper is organized as follows. In section 2, we give the definition of $n$-ary Hom-Nambu algebra and some examples. The definition of derivation is extended to $n$-ary Hom-Nambu algebras from Hom-Lie algebras, and all derivations structure a Hom-Lie algebra in section 3. The most important part is section 4, we define $\omega$-derivations, and study some properties of $\omega$-derivations.
 We can structure a new $n$-ary Hom-Nambu algebra from a $(kn+1)$-ary Hom-Nambu algebra with these properties. In section 5, we also give representations of $n$-ary Hom-Nambu algebras.

\section{$n$-ary Hom-Nambu algebra}
\bdefn $[2]$
A Leibniz $n$-algebra is a pair $(\g, [\cdot,\cdots,\cdot])$ consisting of a vector space $\g$ and a multilinear map $[\cdot,\cdots,\cdot]: \underbrace{\g\times\cdots\times \g}_{n}\rightarrow \g,$ satisfying
\begin{align*}
\begin{split}
[[x_1,\cdots,x_n],y_1,\cdots,y_{n-1}]\\
=\sum_{i=1}^n[x_1,\cdots,x_{i-1},&[x_i,y_1,\cdots, y_{n-1}],x_{i+1},\cdots,x_n].
\end{split}\end{align*}
\edefn

\bdefn $[7]$
An $n$-ary Hom-Nambu algebra is a triple $(\g, [\cdot,\cdots,\cdot],\alpha)$ consisting of a vector space $\g$, a multilinear
 map $[\cdot,\cdots,\cdot]: \underbrace{\g\times\cdots\times \g}_{n}\rightarrow \g$ and a family $\alpha=(\alpha_{i})_{1\leq i\leq n-1}$
  of linear maps $\alpha_{i}: \g\rightarrow \g,$ satisfying
\begin{align}
\begin{split}
[[x_1,\cdots,x_n],\alpha_{1}(y_1),\cdots,\alpha_{n-1}(y_{n-1})]\\
=\sum_{i=1}^n[\alpha_{1}(x_1),\cdots,\alpha_{i-1}(x_{i-1}),&[x_i,y_1,\cdots, y_{n-1}],\alpha_{i}(x_{i+1}),\cdots,\alpha_{n-1}(x_n)].
\end{split}\end{align}
An $n$-ary Hom-Nambu algebra $(\g, [\cdot,\cdots,\cdot],\alpha)$ is multiplicative, if
 $\alpha=(\alpha_{i})_{1\leq i\leq n-1}$\\ with $\alpha_{1}=\cdots=\alpha_{n-1}=\alpha$ and satisfying
$$\alpha([x_{1},\cdots,x_{n}])=[\alpha(x_{1}),\cdots,\alpha(x_{n})], \forall x_{1},x_{2},\cdots,x_{n}\in \g.$$
\edefn
If the $n$-ary Hom-Nambu algebra $(\g, [\cdot,\cdots,\cdot],\alpha)$ is multiplicative, then the equation (2.1) can be read:
\begin{align}
\begin{split}
[[x_1,\cdots,x_n],\alpha(y_1),\cdots,\alpha(y_{n-1})]\\
=\sum_{i=1}^n[\alpha(x_1),\cdots,\alpha(x_{i-1}),&[x_i,y_1,\cdots, y_{n-1}],\alpha(x_{i+1}),\cdots,\alpha(x_n)].
\end{split}\end{align}
\bdefn $[5]$
A Hom-Lie algebra is a triple  $(\g, [\cdot, \cdot]_\g, \alpha)$ consisting of a  vector space $\g$, a bilinear map (bracket) $[\cdot, \cdot]_\g:\g\times \g\rightarrow \g$ and a map $\alpha:\g\rightarrow \g$ satisfying
\begin{align*}[x, y]=-[y, x],\end{align*}
$$[\alpha(x), [y, z]]+[\alpha(y), [z, x]]+[\alpha(z), [x, y]]=0, \qquad\forall x, y, z\in\g.$$
\edefn
\bdefn $[5]$
Let $(\g,[\cdot,\cdot],\alpha)$ be a Hom-Lie algebra. Then $I$ is a Hom-ideal of $\g$, if $I$ satisfies $[I,\g]\subseteq I$ and $\alpha(I)\subseteq I$.
\edefn
\bdefn
Let $(\g, [\cdot,\cdots,\cdot],\alpha)$ and $(\g^{'}, [\cdot,\cdots,\cdot]^{'},\alpha^{'})$ be two $n$-ary Hom-Nambu algebras,
 where $\alpha=(\alpha_{i})_{1\leq i\leq n-1}$ and $\alpha^{'}=(\alpha^{'}_{i})_{1\leq i\leq n-1}.$ A linear map $f: \g\rightarrow \g$ is an
 $n$-ary Hom-Nambu algebra morphism if it satisfies
$$f[x_{1},\cdots,x_{n}]=[f(x_{1}),\cdots,f(x_{n})]^{'},$$
$$f\circ \alpha_{i}=\alpha_{i}^{'}\circ f, \forall i=1,\cdots,n-1.$$
\edefn

\beg
Let $(\g, [\cdot,\cdots,\cdot])$ be a Leibniz $n$-algebra and let $\rho: \g\rightarrow \g$ be a Leibniz $n$-algebra morphism.
Then $(\g, \rho\circ[\cdot,\cdots,\cdot],\rho)$ is a multiplicative $n$-ary Hom-Nambu algebra.
\bpf
Put $[\cdot,\cdots,\cdot]_{\rho}:=\rho\circ[\cdot,\cdots,\cdot].$ Then
\begin{align*}
&\rho[x_{1},\cdots,x_{n}]_{\rho}=\rho(\rho[x_{1},\cdots,x_{n}])\\
=&\rho[\rho(x_{1}),\cdots, \rho(x_{n})]\\
=&[\rho(x_{1}),\cdots, \rho(x_{n})]_{\rho}.
\end{align*}
Moreover, we have
\begin{align*}
&[[x_{1},\cdots,x_{n}]_{\rho},\rho(y_{1}),\cdots,\rho(y_{n-1})]_{\rho}\\
=&\rho[[x_{1},\cdots,x_{n}]_{\rho},\rho(y_{1}),\cdots,\rho(y_{n-1})]\\
=&\rho[\rho[x_{1},\cdots,x_{n-1}],\rho(y_{1}),\cdots,\rho(y_{n-1})]\\
=&\rho^{2}[[x_{1},\cdots,x_{n}],y_{1},\cdots,y_{n-1}]\\
=&\rho^{2}\sum_{i=1}^{n}[x_{1},\cdots,[x_{i},y_{1},\cdots,y_{n-1}],\cdots,x_{n}]\\
=&\sum_{i=1}^{n}[\rho(x_{1}),\cdots,[x_{i},y_{1},\cdots,y_{n-1}]_{\rho},\cdots,\rho(x_{n})]_{\rho}.
\end{align*}
Therefore, $(\g, \rho\circ[\cdot,\cdots,\cdot],\rho)$ is a multiplicative $n$-ary Hom-Nambu algebra.
\epf
\eeg
Similarly, we have
\beg
Let $(\g, [\cdot,\cdots,\cdot],\alpha)$ be a multiplicative $n$-ary Hom-Nambu algebra and let $\beta: \g\rightarrow \g$ be an $n$-ary Hom-Nambu algebra morphism. Then $(\g, \beta\circ[\cdot,\cdots,\cdot],\beta\circ\alpha)$ is a multiplicative $n$-ary Hom-Nambu algebra.
\eeg

\section{Derivations of $n$-ary Hom-Nambu algebra}
Let $(\g, [\cdot,\cdots,\cdot]_g, \alpha)$ be a multiplicative $n$-ary Hom-Nambu algebra. For any nonnegative integer $k$, denote by $\alpha^k$ the $k$-times composition of $\alpha$,  i.e.\\
$$\alpha^k=\alpha\circ\cdots\circ\alpha \quad(k-times).$$
In particular, $\alpha^0=\mathrm{Id}$ and $\alpha^1=\alpha$.
\begin{defn}
For any nonnegative integer $k$, a linear map $D:\g \rightarrow \g$ is called an $\alpha^k$-derivation of the multiplicative  $n$-ary Hom-Nambu algebra $(\g, [\cdot,\cdots,\cdot]_\g, \alpha)$, if
\begin{equation*} [D, \alpha]=0, \quad i.e.\quad{D}\circ\alpha=\alpha\circ{D}, \end{equation*}
and
\begin{equation*}D[x_1,\cdots,x_n ]_\g=\sum_{i=1}^n[\alpha^{k}(x_1),\cdots,D(x_i),\cdots,\alpha^{k}(x_1)]_\g. \end{equation*}
\end{defn}

Denote by $ \texttt{\Der}_{\alpha^s}(\g)$ is the set of $\alpha^s$-derivations of the multiplicative  $n$-ary Hom-Nambu algebra $(\g, [\cdot,\cdots,\cdot]_\g, \alpha)$.
For any $D \in \texttt{\Der}_{\alpha^{k}}(\g)$ and $D^{'} \in \texttt{\Der}_{\alpha^{s}}(\g), $ define their commutator $[D, D^{'}]$ as usual:
\begin{equation*}[D, D^{'}]=D \circ D^{'}-D^{'} \circ D.\end{equation*}
\begin{lem}\label{lemma3.2}For any $D \in \texttt{\Der}_{\alpha^{k}}(\g)$ and $D^{'} \in \texttt{\Der}_{\alpha^{s}}(\g)$, we have\\
$$[D, D^{'}]\in \texttt{\Der}_{\alpha^{k+s}}(\g).$$
\end{lem}
Denote by \begin{eqnarray*}\texttt{\Der}(\g)=\oplus_{k\geq 0} \texttt{\Der}_{\alpha^{k}}(\g).\end{eqnarray*}
By Lemma 3.2,  obviously we have
\begin{prop}With the above notations,  $({\textbf{\Der}(\g)},[\cdot,\cdot],\alpha^{'})$ is a Hom-Lie algebra, with $\alpha^{'}(D)=D\circ\alpha$. \end{prop}

At the end of this section,  we consider the derivation extension of the multiplicative Hom-Leibniz algebra $(\g, [\cdot, \cdot]_\g, \alpha)$ and give an application of the $\alpha$-derivation $ \texttt{\Der}_{\alpha}(\g)$.

For any  linear map $D:\g \rightarrow \g$,  consider the vector space $\g \oplus RD$. Define a multilinear bracket operation$[\cdot,\cdot]_{\g}$ on $\g \oplus RD$ by
$$[x_1+m D,x_2+n D]_{D}=[x_1, x_2]_{\g}+mD(x_2)-{nD(x_1)} \qquad \forall x_1, x_2\in \g.$$
Define a linear map $\alpha^{'}:\g \oplus RD \rightarrow \g \oplus RD$ by $\alpha^{'}(x_1+mD)=\alpha(x_1)+mD, $ i.e.\\
\begin{displaymath}
\mathbf{\alpha^{'}} =
\left( \begin{array}{cc}
\alpha & 0 \\
0 & 1  \\
\end{array}\right).
\end{displaymath}
\begin{prop}
With the above notations, if $(\g \oplus RD,  [\cdot, \cdot]_{D}, \alpha^{'})$ is a multiplicative Hom-Leibniz algebra, then D is an $\alpha$-derivation of the multiplicative Hom-Leibniz algebra $(\g, [\cdot, \cdot]_\g, \alpha)$.
\end{prop}
\bpf
Since $(\g \oplus RD,  [\cdot, \cdot]_{D}, \alpha^{'})$ is a multiplicative Hom-Leibniz algebra, so we have
\begin{align*}
&[[x, y]_{D},\alpha^{'}(D)]_{D}=[\alpha^{'}(x), [y, D]_{D}]_{D}+[[x,D]_{D},\alpha^{'}(y)]_{D},
\end{align*}
that is $D([x,y])=[\alpha(x),D(y)]+[D(x),\alpha(y)].$
\epf
For $\x\in \g^{\wedge^{n-1}}$ satisfying $\alpha(\x)=\x$ and $k\geq 0,$ we define the map $\mathrm{ad}_{k}(\x)\in \mathrm{End}(\g)$ by
$$\mathrm{ad}_{k}(\x)(y)=[\alpha^{k}(y),x_{1},\cdots,x_{n-1}], \forall y\in \g.$$

\blem
The map $\mathrm{ad}_{k}(\x)$ is an $\alpha^{k+1}$-derivation and is called an inner $\alpha^{k+1}$-derivation.
\elem
\bpf For any $y_1,\cdots,y_n\in \g$, we have
\begin{align*}
&\mathrm{ad}_{k-1}(\x)[y_1,\cdots,y_n]\\
=&[\alpha^{k}[y_1,\cdots,y_n],x_1,\cdots,x_{n-1}]\\
=&[[\alpha^{k}(y_1),\cdots,\alpha^{k}(y_n)],\alpha(x_1),\cdots,\alpha(x_{n-1})]\\
=&\sum_{i=1}^{n}[\alpha^{k+1}(y_1),\cdots,\alpha^{k+1}(y_{i-1}),[\alpha^{k}(y_i),\alpha(x_1),\cdots,\alpha(x_{n-1})],\cdots,\alpha^{k+1}(y_n)]\\
=&\sum_{i=1}^{n}[\alpha^{k+1}(y_1),\cdots,\alpha^{k+1}(y_{i-1}),\mathrm{ad}_{k-1}(\x)(y_i),\cdots,\alpha^{k+1}(y_n)].
 \end{align*}
Then $\mathrm{ad}_{k-1}(\x)\in \texttt{\Der}_{\alpha^{k+1}}(\g)$.
\epf
We denote by $\mathrm{Inn}_{\alpha^{k}}(\g)$ the $\K$-vector space generated by all inner $\alpha^{k+1}$-derivations.
\bprop
The $\mathrm{Inn}_{\alpha^{k}}(\g)$ is a Hom-ideal of ${\textbf{\Der}(\g)}$.
\eprop
\bpf $(\mathrm{Der}(\g), [\cdot,\cdot],\alpha^{'})$ is a Hom-Lie algebra. We show that $\mathrm{Inn}(\g)$
is a Hom-ideal. Let $\mathrm{ad}_{k-1}(\x)(y)=[\alpha^{k-1}(y),x_{1},\cdots,x_{n-1}]$ be an inner $\alpha^{k}$-derivation on $\g$ and $D\in \mathrm{Der}_{\alpha^{k^{'}}}(\g)$ for $k\geq 1$ and $k^{'}\geq 0.$ Then
 $[D,\mathrm{ad}_{k-1}(\x)]\in \mathrm{Der}_{\alpha^{k+k^{'}}}(\g)$
and for any $y\in \g$
\begin{align*}
&[D,\mathrm{ad}_{k-1}(\x)](y)\\
=&D[\alpha^{k-1}(y),x_{1},\cdots,x_{n-1}]-[\alpha^{k-1}(D(y)),x_{1},\cdots,x_{n-1}]\\
=&[D\alpha^{k-1}(y),\alpha^{k^{'}}(x_1),\cdots,\alpha^{k^{'}}(x_{n-1})]\\
&+\sum_{i=1}^{n-1}[\alpha^{k+k^{'}-1}(y),\alpha^{k^{'}}(x_1),\cdots,D(x_i),\cdots,\alpha^{k^{'}}(x_{n-1})]\\
&-[\alpha^{k-1}(Dy),\alpha^{k^{'}}(x_1),\cdots,\alpha^{k^{'}}(x_{n-1})]\\
=&\sum_{i=1}^{n-1}[\alpha^{k+k^{'}-1}(y),\alpha^{k^{'}}(x_1),\cdots,D(x_i),\cdots,\alpha^{k^{'}}(x_{n-1})]\\
=&\sum_{i=1}^{n-1}[\alpha^{k+k^{'}-1}(y),x_1,\cdots,D(x_i),\cdots,x_{n-1}]\\
=&\sum_{i=1}^{n-1} \mathrm{ad}_{k+k^{'}-1}(x_{1}\wedge\cdots\wedge D(x_{i})\wedge\cdots\wedge x_{n-1})(y).
\end{align*}
Therefore, $[D,\mathrm{ad}_{k-1}(\x)]\in \mathrm{Inn}_{\alpha^{k+k^{'}-1}}(\g).$
And we have
$$\alpha^{'}(\mathrm{ad}_{k-1}(\x))(y)=\mathrm{ad}_{k-1}(\x)\circ\alpha(y)=\mathrm{ad}_{k}(\x)(y).$$
So the $\mathrm{Inn}_{\alpha^{k}}(\g)$ is a Hom-ideal of ${\textbf{\Der}(\g)}$.
\epf

\section{$\omega$-Derivations}
\bdefn
Let $\A$ be a vector space equipped with an $n$-linear map $\omega:\A^{\otimes n}\rightarrow \A$, $\alpha$ is a linear map on $\A$. A map $f:\A\rightarrow \A$ is called an $\omega$-$\alpha^{k}$-derivation if it satisfies
\begin{align*}
f(\omega(a_1,\cdots,a_n))=\sum_{i=1}^{n}\omega(\alpha^{k}(a_1),\cdots,f(a_i),\cdots,\alpha^{k}(a_n));
\end{align*}
\begin{align*}
f\circ\alpha=\alpha\circ f.
\end{align*}
\edefn
Denote by $ \texttt{\Der}_{\alpha^k}^{\omega}(\A)$ the set of $\omega$-$\alpha^k$-derivations of the vector space $\A$.
For any $f \in \texttt{\Der}_{\alpha^k}^{\omega}(\A)$ and $g \in \texttt{\Der}_{\alpha^s}^{\omega}(\A), $ define their commutator $[f, g]$ as usual:
\begin{equation*}[f,g]=f\circ g-g\circ f.\end{equation*}
\begin{lem}For any $f \in \texttt{\Der}_{\alpha^k}^{\omega}(\A)$ and $g \in \texttt{\Der}_{\alpha^s}^{\omega}(\A), $ we have\\
$$[f,g]\in \texttt{\Der}_{\alpha^{k+s}}^{\omega}(\A).$$
\end{lem}
Denote by \begin{eqnarray*}\texttt{\Der}^{\omega}(\A)=\oplus_{k\geq 0} \texttt{\Der}_{\alpha^{k}}^{\omega}(\A).\end{eqnarray*}
By Lemma 4.2,  obviously we have
\begin{prop}With the above notations,  $({\textbf{\Der}^{\omega}(\A)},[\cdot,\cdot],\alpha^{'})$ is a Hom-Lie algebra, with $\alpha^{'}(f)=f\circ\alpha$. \end{prop}

\begin{lem}
If $f\in \texttt{\Der}_{\alpha^k}^{\omega}(\A)$ and $f\in \textbf{\Der}_{\alpha^k}^{\sigma}(\A)$, then $f\in \textbf{\Der}_{\alpha^k}^{\omega+\sigma}(\A)$. Here $\sigma:\A^{\otimes n}\rightarrow \A$ is also an $n$-linear map.
\end{lem}
\bpf
If $f\in \textbf{\Der}^{\omega}(\A)$ and $f\in \textbf{\Der}^{\sigma}(\A)$, then we have
\begin{align*}
&f((\omega+\delta)(a_1,\cdots,a_n))\\
=&f(\omega(a_1,\cdots,a_n))+f(\sigma(a_1,\cdots,a_n))\\
=&\sum_{i=1}^{n}\omega(\alpha^{k}(a_1),\cdots,f(a_i),\cdots,\alpha^{k}(a_n))+\sum_{i=1}^{n}\sigma(\alpha^{k}(a_1),\cdots,f(a_i),\cdots,\alpha^{k}(a_n))\\
=&\sum_{i=1}^{n}(\omega+\sigma)(\alpha^{k}(a_1),\cdots,f(a_i),\cdots,\alpha^{k}(a_n)).
\end{align*}
So we get the conclusion.
\epf
\beg
Let $(\g, [\cdot,\cdots,\cdot],\alpha)$ be a multiplicative $n$-ary Hom-Nambu algebra. we define $\omega:\g^{\otimes n}\rightarrow \g$
\begin{align*}\begin{split}
\omega(x_1,\cdots,x_n)=[x_1,\cdots,x_n]
 \end{split}\end{align*}
for any $x_1,\cdots,x_n\in \g.$
If $f\in \textbf{\Der}(\g)$, then $f\in \textbf{\Der}^{\omega}(\g)$.
\eeg
\begin{prop}
Let $\omega_{i}:\A^{\otimes n_{i}}\rightarrow \A$ be $n_{i}$-linear maps satisfying $\omega_{i}\circ\alpha=\alpha\circ\omega_{i}$ for $i=1,\cdots,k$ and
 let $\omega:\A^{\otimes k}\rightarrow \A$ be a $k$-linear map. If $f\in \bigcap_{i=1}^{k}\textbf{\Der}_{\alpha^{t}}^{\omega_{i}}(\A)\bigcap\textbf{\Der}_{\alpha^{t}}^{\omega}(\A)$, then $f\in \textbf{\Der}_{\alpha^{t}}^{\sigma}(\A)$ with $\sigma:\A^{\otimes n}\rightarrow \A$. Here $n=n_1+\cdots+n_k$, $\sigma(a_1,\cdots,a_n)\triangleq\omega(\omega_{1}(a_1,\cdots,a_{n_{1}}),\cdots,\omega_{k}(a_1,\cdots,a_{n_{k}}))$ and $s=n-n_k+1=n_1+\cdots+n_{k-1}+1$.
\end{prop}\bpf
On the one hand, we have
\begin{align*}
&f(\sigma(a_1,\cdots,a_n))\\
=&f(\omega(\omega_{1} (a_1,\cdots,a_{n_1}),\cdots,\omega_{k}(a_s,\cdots,a_n)))\\
=&\sum_{i=1}^{k}\omega(\alpha^{t}(\omega_1(a_1,\cdots,a_{n_1})),\cdots,
f(\omega_i(a_{m_1},\cdots,a_{m_{n_i}})),\cdots,\alpha^{t}(\omega_k(a_s,\cdots,a_{n})))\\
=&\sum_{i=1}^{k}\omega(\omega_1(\alpha^{t}(a_1),\cdots,\alpha^{t}(a_{n_1}))),\cdots,
f(\omega_i(a_{m_1},\cdots,a_{m_{n_i}})),\cdots,\omega_k(\alpha^{t}(a_s),\cdots,\alpha^{t}(a_{n}))),
\end{align*}
and here $$f(\omega_i(a_{m_1},\cdots,a_{m_{n_i}}))=\sum_{j=1}^{m_{n_i}}\omega_i(\alpha^{t}(a_{m_1}),\cdots,f(a_{m_{n_j}}),
\cdots,\alpha^{t}(a_{m_{n_i}})).$$
On the other hand,
\begin{align*}
&\sum_{i=1}^{n}\sigma(\alpha^{t}(a_1),\cdots,f(a_i),\cdots,\alpha^{t}(a_n))\\
=&\sum_{i=1}^{n}\omega(\omega_1(\alpha^{t}(a_1),\cdots,\alpha^{t}(a_{n_1}))),\cdot,\omega_i
(\alpha^{t}(a_{p_1}),\cdots,f(a_{p_t}),\cdot),\cdots,\omega_k(\alpha^{t}(a_s),\cdots,\alpha^{t}(a_{n}))).
\end{align*}
So we have
$f(\sigma(a_1,\cdots,a_n))=\sum_{i=1}^{n}\sigma(\alpha^{t}(a_1),\cdots,f(a_i),\cdots,\alpha^{t}(a_n)).$
\epf

\begin{lem}
Let $\omega:\A^{\otimes (n+1)}\rightarrow \A$ be an $(n+1)$-linear map and let $\mu_{i}:\g\times \g\rightarrow \g$ be bilinear maps given by
\begin{align*}
&\mu_{i}(a_{1}\otimes\cdots\otimes a_{n},b_{1}\otimes\cdots\otimes b_{n})
=\alpha^{k}(a_1)\otimes\cdots\otimes\omega(a_i,b_1,\cdots,b_n)\otimes\cdots\otimes\alpha^{k}(a_n).
\end{align*}
Here $\g=\A^{\otimes n}$ and $1\leq i\leq n$. Suppose that $f\in \textbf{\Der}_{\alpha^{0}}^{\omega}(\A)$ and $\varphi:\g\rightarrow \g$ is given by
\begin{align*}\begin{split}
\varphi(a_{1}\otimes\cdots\otimes a_{n})=\sum_{j=1}^{n}a_{1}\otimes\cdots\otimes f(a_i)\otimes\cdots\otimes a_{n}.
\end{split}\end{align*}
Then $\varphi\in\textbf{\Der}_{\alpha^{0}}^{\mu_{i}}(\A)$.
\end{lem}
\bpf
On the one hand, we have
\begin{align*}
&\varphi(\mu_{i}(a_{1}\otimes\cdots\otimes a_{n},b_{1}\otimes\cdots\otimes b_{n}))\\
=&\varphi(\alpha^{k}(a_1)\otimes\cdots\otimes\omega(a_i,b_1,\cdots,b_n)\otimes\cdots\otimes\alpha^{k}(a_n))\\
=&\sum_{j=1}^{i-1}\alpha^{k}(a_1)\otimes\cdots\otimes f\alpha^{k}(a_j)\otimes\cdots\otimes\omega(a_i,b_1,\cdots,b_n)\otimes\cdots\otimes\alpha^{k}(a_n)\\
+&\sum_{j=i+1}^{n}\alpha^{k}(a_1)\otimes\cdots\otimes\omega(a_i,b_1,\cdots,b_n)\otimes\cdots\otimes f\alpha^{k}(a_j)\otimes\cdots\otimes\alpha^{k}(a_n)\\
+&\alpha^{k}(a_1)\otimes\cdots\otimes f(\omega(a_i,b_1,\cdots,b_n))\otimes\cdots\otimes\alpha^{k}(a_n)\\
=&\sum_{j=1}^{i-1}\alpha^{k}(a_1)\otimes\cdots\otimes \alpha^{k}f(a_j)\otimes\cdots\otimes\omega(a_i,b_1,\cdots,b_n)\otimes\cdots\otimes\alpha^{k}(a_n)\\
+&\sum_{j=i+1}^{n}\alpha^{k}(a_1)\otimes\cdots\otimes\omega(a_i,b_1,\cdots,b_n)\otimes\cdots\otimes \alpha^{k}f(a_j)\otimes\cdots\otimes\alpha^{k}(a_n)\\
+&\sum_{j=1}^{n}\alpha^{k}(a_1)\otimes\cdots\otimes\omega(a_i,b_1,\cdots,f(b_j),\cdots,b_n)\otimes\cdots\otimes\alpha^{k}(a_n)\\
+&\alpha^{k}(a_1)\otimes\cdots\otimes\omega(f(a_i),b_1,\cdots,b_n)\otimes\cdots\otimes\alpha^{k}(a_n).
\end{align*}
The other hand,
\begin{align*}\begin{split}
&\mu_{i}(\varphi(a_{1}\otimes\cdots\otimes a_{n}),b_{1}\otimes\cdots\otimes b_{n})+\mu_{i}(a_{1}\otimes\cdots\otimes a_{n},\varphi(b_{1}\otimes\cdots\otimes b_{n}))\\
=&\mu_{i}(\sum_{j=1}^{n}a_{1}\otimes\cdots\otimes f(a_j)\otimes\cdots\otimes a_{n},b_{1}\otimes\cdots\otimes b_{n})\\
+&\mu_{i}(a_{1}\otimes\cdots\otimes a_{n},\sum_{j=1}^{n}b_{1}\otimes\cdots\otimes f(b_j)\otimes\cdots\otimes b_{n})\\
=&\sum_{j=1}^{n}\mu_{i}(a_{1}\otimes\cdots\otimes f(a_j)\otimes\cdots\otimes a_{n},b_{1}\otimes\cdots\otimes b_{n})\\
+&\sum_{j=1}^{n}\mu_{i}(a_{1}\otimes\cdots\otimes a_{n},b_{1}\otimes\cdots\otimes f(b_j)\otimes\cdots\otimes b_{n})\\
=&\sum_{j=1}^{i-1}\alpha^{k}(a_1)\otimes\cdots\otimes \alpha^{k}f(a_j)\otimes\cdots\otimes\omega(a_i,b_1,\cdots,b_n)\otimes\cdots\otimes\alpha^{k}(a_n)\\
+&\sum_{j=i+1}^{n}\alpha^{k}(a_1)\otimes\cdots\otimes\omega(a_i,b_1,\cdots,b_n)\otimes\cdots\otimes \alpha^{k}f(a_j)\otimes\cdots\otimes\alpha^{k}(a_n)\\
+&\alpha^{k}(a_1)\otimes\cdots\otimes\omega(f(a_i),b_1,\cdots,b_n)\otimes\cdots\otimes\alpha^{k}(a_n)\\
+&\sum_{j=1}^{n}\alpha^{k}(a_1)\otimes\cdots\otimes\omega(a_i,b_1,\cdots,f(b_j),\cdots,b_n)\otimes\cdots\otimes\alpha^{k}(a_n).
\end{split}\end{align*}
So we have
\begin{align*}\begin{split}&\varphi(\mu_{i}(a_{1}\otimes\cdots\otimes a_{n},b_{1}\otimes\cdots\otimes b_{n}))\\
=&\mu_{i}(\varphi(a_{1}\otimes\cdots\otimes a_{n}),b_{1}\otimes\cdots\otimes b_{n})
+\mu_{i}(a_{1}\otimes\cdots\otimes a_{n},\varphi(b_{1}\otimes\cdots\otimes b_{n})),\end{split}\end{align*}
that is $\varphi\in\textbf{\Der}_{\alpha^{0}}^{\mu_{i}}(\A)$.
\epf
\begin{thm}
(1) Let $(\g,[\cdot,\cdots,\cdot],\alpha)$ be a multiplicative $(n+1)$-ary Hom-Nambu algebra. Then $D_{n}(\g)$ is a Leibniz algebra with respect to the bracket
$$[a_{1}\otimes\cdots\otimes a_{n},b_{1}\otimes\cdots\otimes b_{n}]\triangleq\sum_{i=1}^{n}a_1
\otimes\cdots\otimes[a_i,b_1,\cdots,b_n]\otimes\cdots\otimes a_n.$$
(2) Moreover, if we define
$$[a_{1}\otimes\cdots\otimes a_{n},b_{1}\otimes\cdots\otimes b_{n}]\triangleq\sum_{i=1}^{n}\alpha(a_1)
\otimes\cdots\otimes[a_i,b_1,\cdots,b_n]\otimes\cdots\otimes\alpha(a_n)$$ and
$\alpha^{'}:D_{n}(\g)\rightarrow D_{n}(\g)$ satisfying
$$\alpha^{'}(a_{1}\otimes\cdots\otimes a_{n})=\alpha(a_{1})\otimes\cdots\otimes\alpha(a_{n}).$$
Then $(D_{n}(\g),[\cdot,\cdot],\alpha^{'})$ is a multiplicative Hom-Leibniz algebra.
\end{thm}
\bpf
(1) From Lemma 4.7, let $k=1$, and
\begin{align*}
&\mu_{i}(a_{1}\otimes\cdots\otimes a_{n},b_{1}\otimes\cdots\otimes b_{n})
=[a_1,\cdots,[a_i,b_1,\cdots,b_n],\cdots,a_n].
\end{align*}
Fix $x_1,\cdots,x_n\in\g,$ define $f(a)=[a, x_1,\cdots,x_n]$ , $\varphi:D_{n}(\g)\rightarrow D_{n}(\g)$ is given by
\begin{align*}
\varphi(a_{1}\otimes\cdots\otimes a_{n})=\sum_{j=1}^{n}a_{1}\otimes\cdots\otimes f(a_i)\otimes\cdots\otimes a_{n}.
\end{align*}
From Lemma 4.7, we get $\varphi\in\textbf{\Der}_{\alpha^{0}}^{\mu_{i}}(\g)$.
And from Lemma 4.4, we know\\ $\varphi\in\textbf{\Der}_{\alpha^{0}}^{\mu}(\g)$.
Here $\mu=\mu_{1}+\cdots+\mu_{n}$, and $$\mu(a_{1}\otimes\cdots\otimes a_{n},b_{1}\otimes\cdots\otimes b_{n})=[a_{1}\otimes\cdots\otimes a_{n},b_{1}\otimes\cdots\otimes b_{n}].$$
$\varphi\in\textbf{\Der}_{\alpha^{0}}^{\mu}(\g)$, that is
\begin{align*}
&\varphi(\mu(a_{1}\otimes\cdots\otimes a_{n},b_{1}\otimes\cdots\otimes b_{n}))=\mu(\varphi(a_{1}\otimes\cdots\otimes a_{n}),b_{1}\otimes\cdots\otimes b_{n})\\
+&\mu(a_{1}\otimes\cdots\otimes a_{n},\varphi(b_{1}\otimes\cdots\otimes b_{n})),
\end{align*}
that is
\begin{align*}
&\varphi([a_{1}\otimes\cdots\otimes a_{n},b_{1}\otimes\cdots\otimes b_{n}])=[\varphi(a_{1}\otimes\cdots\otimes a_{n}),b_{1}\otimes\cdots\otimes b_{n}]\\
+&[a_{1}\otimes\cdots\otimes a_{n},\varphi(b_{1}\otimes\cdots\otimes b_{n})]).
\end{align*}
 It is equivalent to
\begin{align*}
&[[a_{1}\otimes\cdots\otimes a_{n},b_{1}\otimes\cdots\otimes b_{n}],c_{1}\otimes\cdots\otimes c_{n}]\\
=&[[a_{1}\otimes\cdots\otimes a_{n},c_{1}\otimes\cdots\otimes c_{n}],b_{1}\otimes\cdots\otimes b_{n}]\\
+&[a_{1}\otimes\cdots\otimes a_{n},[b_{1}\otimes\cdots\otimes b_{n},c_{1}\otimes\cdots\otimes c_{n}]].
\end{align*}
So $D_{n}(\g)$ is a Leibniz algebra.\\
(2) According to the definition of bracket, we have
\begin{align*}
&[[a_{1}\otimes\cdots\otimes a_{n},b_{1}\otimes\cdots\otimes b_{n}],\alpha(c_{1})\otimes\cdots\otimes\alpha(c_{n})]\\
=&[\sum_{i=1}^{n}\alpha(a_{1})\otimes\cdots\otimes[a_i,b_{1}\otimes\cdots\otimes b_{n}]\otimes\cdots\otimes\alpha(a_{n}),\alpha(c_{1})\otimes\cdots\otimes\alpha(c_{n})]\\
=&\sum_{i=1}^{n}[\alpha(a_{1})\otimes\cdots\otimes[a_i,b_{1}\otimes\cdots\otimes b_{n}]\otimes\cdots\otimes\alpha(a_{n}),\alpha(c_{1})\otimes\cdots\otimes\alpha(c_{n})]\\
=&\sum_{i=1}^{n}\sum_{j=1}^{i-1}\alpha^{2}(a_{1})\otimes\cdots\otimes[\alpha(a_j),\alpha(c_1),\cdots,\alpha(c_n)]
\otimes\cdots\otimes\alpha([a_i,b_{1},\\
&\cdots, b_{n}])\otimes\cdots\otimes\alpha^{2}(a_{n})\tag{${1}$}\\
+&\sum_{i=1}^{n}\sum_{j=i+1}^{n}\alpha^{2}(a_{1})\otimes\cdots\otimes\alpha([a_i,b_{1},\cdots, b_{n}])
\otimes\cdots\otimes[\alpha(a_j),\alpha(c_1),\\
&\cdots,\alpha(c_n)]\otimes\cdots\otimes\alpha^{2}(a_{n})\tag{${2}$}\\
+&\sum_{i=1}^{n}\alpha^{2}(a_{1})\otimes\cdots\otimes[[a_i,b_{1}\otimes\cdots\otimes b_{n}],\alpha(c_1),\cdots,\alpha(c_n)]\otimes\cdots\otimes\alpha^{2}(a_{n}).\tag{${3}$}
\end{align*}
In the same way,
\begin{align*}
&[[a_{1}\otimes\cdots\otimes a_{n},c_{1}\otimes\cdots\otimes c_{n}],\alpha(b_{1})\otimes\cdots\otimes\alpha(b_{n})]\\
=&\sum_{i=1}^{n}\sum_{j=1}^{i-1}\alpha^{2}(a_{1})\otimes\cdots\otimes[\alpha(a_j),\alpha(b_1),\cdots,\alpha(b_n)]
\otimes\cdots\otimes\alpha([a_i,c_{1}\otimes\cdots\otimes c_{n}])\otimes\\
&\cdots\otimes\alpha^{2}(a_{n})\tag{${2}^{'}$}\\
+&\sum_{i=1}^{n}\sum_{j=i+1}^{n}\alpha^{2}(a_{1})\otimes\cdots\otimes\alpha([a_i,c_{1}\otimes\cdots\otimes c_{n}])
\otimes\cdots\otimes[\alpha(a_j),\alpha(b_1),\cdots,\alpha(b_n)]\otimes\\
&\cdots\otimes\alpha^{2}(a_{n})\tag{${1}^{'}$}\\
+&\sum_{i=1}^{n}\alpha^{2}(a_{1})\otimes\cdots\otimes[[a_i,c_{1}\otimes\cdots\otimes c_{n}],\alpha(b_1),\cdots,\alpha(b_n)]\otimes\cdots\otimes\alpha^{2}(a_{n}).\tag{${3}^{'}$}
\end{align*}
And
\begin{align*}
&[\alpha(a_{1})\otimes\cdots\otimes\alpha(a_{n}),[b_{1}\otimes\cdots\otimes b_{n},c_{1}\otimes\cdots\otimes c_{n}]]\\
=&\sum_{i=1}^{n}[\alpha(a_{1})\otimes\cdots\otimes\alpha(a_{n}),\alpha(b_{1})\otimes\cdots\otimes[b_i,c_1,\cdots,c_n]
\otimes\cdots\otimes \alpha(b_{n})]\\
=&\sum_{i=1}^{n}\sum_{j=1}^{n}\alpha^{2}(a_{1})\otimes\cdots\otimes[\alpha(a_{j}),\alpha(b_{1}),\cdots,[b_i,c_1,\cdots,c_n]
,\cdots,\alpha(b_{n})]\otimes\cdots\otimes\alpha^{2}(a_{n}).\tag{${3}^{''}$}
\end{align*}
Obviously, $(1)=(1^{'}),(2)=(2^{'}),(3)=(3^{'})+(3^{''})$, so \begin{align*}
&[[a_{1}\otimes\cdots\otimes a_{n},b_{1}\otimes\cdots\otimes b_{n}],\alpha(c_{1})\otimes\cdots\otimes\alpha(c_{n})]\\
=&[[a_{1}\otimes\cdots\otimes a_{n},c_{1}\otimes\cdots\otimes c_{n}],\alpha(b_{1})\otimes\cdots\otimes\alpha(b_{n})]\\
+&[\alpha(a_{1})\otimes\cdots\otimes\alpha(a_{n}),[b_{1}\otimes\cdots\otimes b_{n},c_{1}\otimes\cdots\otimes c_{n}]].
\end{align*}
 So $(D_{n}(\g),[\cdot,\cdot],\alpha^{'})$ is a multiplicative Hom-Leibniz algebra.
\epf
Moreover, we can prove
\begin{cor}
 If $\g$ is a $(kn+1)$-ary Hom-Nambu algebra, then $(\g^{\otimes k},[\cdot,\cdots,\cdot],\alpha^{'})$ is a Hom-Leibniz $(n+1)$-algebra
 with respect to the following bracket
\begin{align*}
&[x_{01}\otimes\cdots\otimes x_{0k},\cdots,x_{n1}\otimes\cdots\otimes x_{nk}]\\
&\triangleq[x_{01},\cdots,x_{11},\cdots,x_{1k},\cdots,x_{n1},\cdots,x_{nk}]\otimes\alpha(x_{02})\otimes\cdots\otimes\alpha(x_{0k})\\
+&\cdots+\alpha(x_{01})\otimes\cdots\otimes\alpha(x_{ok-1})\otimes[x_{0k},x_{11},\cdots,x_{nk}]
\end{align*}
and
$\alpha^{'}:\g^{\otimes k}\rightarrow \g^{\otimes k}$ satisfying
$$\alpha^{'}(a_{1}\otimes\cdots\otimes a_{k})=\alpha(a_{1})\otimes\cdots\otimes\alpha(a_{k}).$$
\end{cor}

\section{Representations of $n$-ary Hom-Nambu algebra}
\begin{defn}
Let $(\g,[\cdot,\cdots,\cdot],\alpha)$ be an $n$-ary Hom-Nambu algebra. We say $M$ is a representation of $(\g,[\cdot,\cdots,\cdot],\alpha)$, if $M$ is a
vector space with a multilinear map $[\cdot,\cdots,\cdot]:\g^{\otimes^{i}}\otimes M\otimes\g^{\otimes^{n-1-i}}\longrightarrow M$ satisfying
\begin{align*}
[[x_1,\cdots,x_n],\alpha^{'}(y_1),\cdots,\alpha^{'}(y_{n-1})]\\
=\sum_{i=1}^n[\alpha^{'}(x_1),\cdots,\alpha^{'}(x_{i-1}),&[x_i,y_1,\cdots, y_{n-1}],\alpha^{'}(x_{i+1}),\cdots,\alpha^{'}(x_n)],
\end{align*}
and one of $x_{1},\cdots,x_{n},y_{1},\cdots,y_{n}$ in $M$, others in $\g$.
Here $\alpha^{'}:M\oplus\g\longrightarrow M\oplus\g$ satisfying
$$\alpha^{'}(x)=\alpha(x),\forall x\in\g;$$
$$\alpha^{'}(m)=m,\forall m\in M.$$
\end{defn}
\begin{prop}
Let $M$ be a representation of
$(\g,[\cdot,\cdots,\cdot],\alpha^{'})$, and $f:\g^{\otimes
n}\longrightarrow M$ is an $n$-linear map. We define an $n$-bracket
on $H=M\oplus\g$ by
$$[(m_1,x_1),\cdots,(m_n,x_n)]\triangleq(\sum_{i=1}^{n}[x_1,\cdots,m_i,\cdots,x_n]+f(x_1,\cdots,x_n),[x_1,\cdots,x_n])$$
and $\alpha^{'}:M\oplus\g\longrightarrow M\oplus\g$ satisfying
$$\alpha^{'}(m,x)=(\alpha(x),m),\forall x\in\g, m\in M.$$
Then $(H,[\cdot,\cdots,\cdot],\alpha^{'})$ is an $n$-ary Hom-Nambu algebra if and only if
\begin{align}
\begin{split}
&f([x_1,\cdots,x_n],\alpha(y_1),\cdots,\alpha(y_{n-1}))+[f(x_1,\cdots,x_n),\alpha(y_1),\cdots,\alpha(y_{n-1})]\\
=&\sum_{i=1}^{n}(f(\alpha(x_1),\cdots,[x_i,y_1,\cdots,y_{n-1}],\cdots,\alpha(x_n))\\
+&[\alpha(x_1),\cdots,f([x_i,y_1,\cdots,y_{n-1}]),\cdots,\alpha(x_n)]).
\end{split}
\end{align}
\end{prop}
\bpf
On the one hand, we have
\begin{align*}
&[[(m_1,x_1),\cdots,(m_n,x_n)],\alpha^{'}(p_1,y_1),\cdots,\alpha^{'}(p_{n-1},y_{n-1})]\\
=&[(\sum_{i=1}^{n}[x_1,\cdots,m_i,\cdots,x_n]+f(x_1,\cdots,x_n),[x_1,\cdots,x_n]),\alpha^{'}(p_1,y_1),\cdots,\alpha^{'}(p_{n-1},y_{n-1})]\\
=&(\sum_{j=1}^{n-1}[[x_1,\cdots,x_n],\cdots,p_j,\cdots,\alpha(y_{n-1})]+[\sum_{i=1}^{n}[x_1,\cdots,m_i,\cdots,x_n]+f(x_1,\cdots,x_n),\\
&\alpha(y_1),\cdots,\alpha(y_{n-1})]+f([x_1,\cdots,x_n],\alpha(y_1),\cdots,\alpha(y_{n-1})),[[x_1,\cdots,x_n],\alpha(y_1),\cdots,\alpha(y_{n-1})])\\
=&(\sum_{j=1}^{n-1}[[x_1,\cdots,x_n],\cdots,p_j,\cdots,\alpha(y_{n-1})],\tag{${1}$}\\
&[[x_1,\cdots,x_n],\alpha(y_1),\cdots,\alpha(y_{n-1})])\tag{${2}$}\\
+&(\sum_{i=1}^{n}[[x_1,\cdots,m_i,\cdots,x_n],\alpha(y_1),\cdots,\alpha(y_{n-1})],\tag{${3}$}\\
&[[x_1,\cdots,x_n],\alpha(y_1),\cdots,\alpha(y_{n-1})])\\
+&([f(x_1,\cdots,x_n),\alpha(y_1),\cdots,\alpha(y_{n-1})],[[x_1,\cdots,x_n],\alpha(y_1),\cdots,\alpha(y_{n-1})])\\
+&(f([x_1,\cdots,x_n],\alpha(y_1),\cdots,\alpha(y_{n-1})),[[x_1,\cdots,x_n],\alpha(y_1),\cdots,\alpha(y_{n-1})]).
\end{align*}
The other hand,
\begin{align*}
&\sum_{i=1}^{n}[\alpha^{'}(m_1,x_1),\cdots,\alpha^{'}(m_{i-1},x_{i-1}),[(m_i,x_i),(p_1,y_1),\cdots,(p_{n-1},y_{n-1})],\cdots,\alpha^{'}(m_n,x_n)]\\
=&\sum_{i=1}^{n}[(m_1,\alpha(x_1)),\cdots,(m_{i-1},\alpha(x_{i-1})),[(m_i,x_i),(p_1,y_1),\cdots,(p_{n-1},y_{n-1})],\cdots,(m_n,\alpha(x_n))]\\
=&\sum_{i=1}^{n}[(m_1,\alpha(x_1)),\cdots,(m_{i-1},\alpha(x_{i-1})),(\sum_{j=1}^{n-1}[x_i,\cdots,p_j,\cdots,y_{n-1}]+[m_i,y_1,\cdots,y_{n-1}]+\\ &f(x_i,y_1,\cdots,y_{n-1}),[x_i,y_1,\cdots,y_{n-1}]),\cdots,(m_n,\alpha(x_n))]\\
=&(\sum_{i=1}^{n}\sum_{j=1}^{n-1}[\alpha(x_1),\cdots,\alpha(x_{i-1}),[x_i,\cdots,p_j,\cdots,y_{n-1}],\cdots,\alpha(x_n)],\tag{${1}^{'}$}\\
&\sum_{i=1}^{n}[\alpha(x_1),\cdots,\alpha(x_{i-1}),[x_i,y_1,\cdots, y_{n-1}],\alpha(x_{i+1}),\cdots,\alpha(x_n)])\tag{${2}^{'}$}\\
+&(\sum_{i=1}^{n}[\alpha(x_1),\cdots,\alpha(x_{i-1}),[m_i,\cdots,y_1,\cdots,y_{n-1}],\cdots,\alpha(x_n)],\tag{${3}^{'}$}\\
&\sum_{i=1}^{n}[\alpha(x_1),\cdots,\alpha(x_{i-1}),[x_i,y_1,\cdots, y_{n-1}],\alpha(x_{i+1}),\cdots,\alpha(x_n)])\\
+&(\sum_{i=1}^{n}[\alpha(x_1),\cdots,\alpha(x_{i-1}),f(x_i,y_1,\cdots,y_{n-1}),\alpha(x_{i+1}),\cdots,\alpha(x_n)],\\
&\sum_{i=1}^{n}[\alpha(x_1),\cdots,\alpha(x_{i-1}),[x_i,y_1,\cdots, y_{n-1}],\alpha(x_{i+1}),\cdots,\alpha(x_n)])\\
+&(\sum_{i=1}^{n}f(\alpha(x_1),\cdots,\alpha(x_{i-1}),[x_i,y_1,\cdots,y_{n-1}],\alpha(x_{i+1}),\cdots,\alpha(x_n)),\\
&\sum_{i=1}^{n}[\alpha(x_1),\cdots,\alpha(x_{i-1}),[x_i,y_1,\cdots, y_{n-1}],\alpha(x_{i+1}),\cdots,\alpha(x_n)]).
\end{align*}
Obviously, $(1)=(1^{'}),(2)=(2^{'}),(3)=(3^{'})$, so $(H,[\cdot,\cdots,\cdot],\alpha^{'})$ is
an $n$-ary Hom-Nambu algebra if and only if
\begin{align*} &f([x_1,\cdots,x_n],\alpha(y_1),\cdots,\alpha(y_{n-1}))+[f(x_1,\cdots,x_n),\alpha(y_1),\cdots,\alpha(y_{n-1})]\\
=&\sum_{i=1}^{n}(f(\alpha(x_1),\cdots,[x_i,y_1,\cdots,y_{n-1}],\cdots,\alpha(x_n))\\
+&[\alpha(x_1),\cdots,f([x_i,y_1,\cdots,y_{n-1}]),\cdots,\alpha(x_n)]). \end{align*}
Then we get the conclusion.
\epf
Obviously, if the condition of the proposition 5.2 holds, then we obtain an extension
$$ 0\rightarrow M\rightarrow H\rightarrow \g\rightarrow 0$$ of $n$-ary Hom-Nambu algebras. Moreover,
\begin{prop}
This extension is split in the category of $n$-ary Hom-Nambu algebras if and only if there exists a linear map $h:\g\rightarrow M$ such that
\begin{align*}
\g\circ\alpha^{'}=\alpha^{'}\circ\g;
\end{align*}
\begin{align}
\begin{split}
f(x_1,\cdots,x_n)=\sum_{i=1}^{n}[x_1,\cdots,,\cdots,x_n]-h([x_1,\cdots,x_n]).
\end{split}
\end{align}
\end{prop}
An easy consequence of these facts is the following natural bejection:
\begin{align}
\begin{split}
Ext(\g,M)\cong Z(\g,M)/B(\g,M).
\end{split}
\end{align}
Here $Ext(\g,M)$ is the set of isomorphism classes of extensions of $\g$ by $M$, $Z(\g,M)$ is the set of all linear maps
$f:\g^{\otimes n}\longrightarrow M$ satisfying (5.1), and $Z(\g,M)$ is the set of such $f$ which satisfies (5.2) for some $k$-linear map $h:\g\rightarrow M$.

We just to prove $B(\g,M)\subseteq Z(\g,M)$. In fact, if $f$ satisfies (5.2), then
\begin{align*}
&f([x_1,\cdots,x_n],\alpha(y_1),\cdots,\alpha(y_{n-1}))+[f(x_1,\cdots,x_n),\alpha(y_1),\cdots,\alpha(y_{n-1})]\\
=&[h([x_1,\cdots,x_n]),\alpha(y_1),\cdots,\alpha(y_{n-1})]\\
+&\sum_{i=1}^{n-1}[[x_1,\cdots,x_n],\alpha(y_1),\cdots,h(\alpha(y_i)),\cdots,\alpha(y_{n-1})]-h([[x_1,\cdots,x_n],\alpha(y_1),\cdots,\alpha(y_{n-1})])\\
+&\sum_{i=1}^{n}[[x_1,\cdots,h(x_i),\cdots,x_n],\alpha(y_1),\cdots,\alpha(y_{n-1})]-[h([x_1,\cdots,x_n]),\alpha(y_1),\cdots,\alpha(y_{n-1})]\\
=&\sum_{i=1}^{n-1}[[x_1,\cdots,x_n],\alpha(y_1),\cdots,h(\alpha(y_i)),\cdots,\alpha(y_{n-1})]\tag{${1}$}\\
-&h([[x_1,\cdots,x_n],\alpha(y_1),\cdots,\alpha(y_{n-1})])\tag{${2}$}\\
+&\sum_{i=1}^{n}[[x_1,\cdots,h(x_i),\cdots,x_n],\alpha(y_1),\cdots,\alpha(y_{n-1})].\tag{${3}$}
\end{align*}
And
\begin{align*}
&\sum_{i=1}^{n}(f(\alpha(x_1),\cdots,[x_i,y_1,\cdots,y_{n-1}],\cdots,\alpha(x_n))\\
+&[\alpha(x_1),\cdots,f([x_i,y_1,\cdots,y_{n-1}]),\cdots,\alpha(x_n)])\\
=&\sum_{i=1}^{n}\sum_{j=1}^{i-1}[\alpha(x_1),\cdots,h(\alpha(x_j)),\cdots,[x_i,y_1,\cdots,y_{n-1}],\cdots,\alpha(x_n)]\tag{${3}^{'}$}\\
+&\sum_{i=1}^{n}\sum_{j=i+1}^{n}[\alpha(x_1),\cdots,[x_i,y_1,\cdots,y_{n-1}],\cdots,h(\alpha(x_j)),\cdots,\alpha(x_n)]\tag{${3}^{''}$}\\
+&\sum_{i=1}^{n}[\alpha(x_1),\cdots,h([x_i,y_1,\cdots,y_{n-1}]),\cdots,\alpha(x_n)]\tag{${4}$}\\
-&\sum_{i=1}^{n}h([\alpha(x_1),\cdots,[x_i,y_1,\cdots,y_{n-1}],\cdots,\alpha(x_n)])\tag{${2}^{'}$}\\
+&\sum_{i=1}^{n}\sum_{j=1}^{n-1}[\alpha(x_1),\cdots,\alpha(x_{i-1}),[x_i,\cdots,h(y_j),\cdots,y_{n-1}],\alpha(x_{i+1}),\cdots,\alpha(x_n)]\tag{${1}^{'}$}\\
+&\sum_{i=1}^{n}[\alpha(x_1),\cdots,\alpha(x_{i-1}),[h(x_i),\cdots,y_{n-1}],\alpha(x_{i+1}),\cdots,\alpha(x_n)]\tag{${3}^{''''}$}\\
-&\sum_{i=1}^{n}[\alpha(x_1),\cdots,\alpha(x_{i-1}),h([x_i,\cdots,y_{n-1}]),\alpha(x_{i+1}),\cdots,\alpha(x_n)].\tag{${4}^{'}$}
\end{align*}
Since $(1)=(1^{'}),(2)=(2^{'}),(3)=(3^{'})+(3^{''})+(3^{'''}),(4)=-(4^{'})$, so $f$ satisfies (5.1).
\beg
Let $(\g,[\cdot,\cdots,\cdot],\alpha^{'})$ be an $(n+1)$-ary Hom-Nambu algebra, M is a representation of $\g$. $Hom(\g,M)$ is the set of isomorphism of $\g\rightarrow M$. One defines the maps
\begin{align*}
&[-,-]:Hom(\g,M)\otimes\rightarrow Hom(\g,M)\\
&[-,-]:D_n(\g)\otimes Hom(\g,M)\rightarrow Hom(\g,M)
\end{align*}
by
\begin{align*}
&[f,x_{1}\otimes\cdots\otimes x_{n}](x)=-[f(x),x_{1}\otimes\cdots\otimes x_{n}],\\
&[x_{1}\otimes\cdots\otimes x_{n},f](x)=-[f(x),x_{1}\otimes\cdots\otimes x_{n}].
\end{align*}
Then $(Hom(\g,M),[-,-])$ is a representation of $D_n(\g)$, here $D_n(\g)=\g^{\otimes n}$.
\eeg

\end{document}